\documentclass[11pt, twoside, a4paper]{article}
\usepackage[english]{babel}
\usepackage{amsmath,amssymb,amsthm, amscd,eucal,epsfig}
\begin{document}

\newtheorem{theorem}{Theorem}
\newtheorem{definition}[theorem]{Definition}
\newtheorem{proposition}[theorem]{Proposition}
\newtheorem{lemma}[theorem]{Lemma}
\newtheorem{corollary}[theorem]{Corollary}
\newtheorem{remark}[theorem]{Remark}
\newtheorem{notation}[theorem]{Notation}
\newtheorem{application}[theorem]{Application}
\newtheorem{definition-proposition}[theorem]{Definition-Proposition}

\title{On calibrated and separating sub-actions}
\author{
Eduardo Garibaldi\\
\footnotesize{Departamento de Matem\'atica}\\
\footnotesize{Universidade Estadual de Campinas}\\
\footnotesize{13083-859 Campinas -- SP, Brasil}\\
\footnotesize{\texttt{garibaldi@ime.unicamp.br}}
\and
Artur O. Lopes\thanks{Partially supported by CNPq, PRONEX -- Sistemas Din\^amicos, Instituto do Mil\^enio, and beneficiary of CAPES financial support.}\\
\footnotesize{Instituto de Matem\'atica}\\
\footnotesize{Universidade Federal do Rio Grande do Sul}\\
\footnotesize{91509-900 Porto Alegre -- RS, Brasil}\\
\footnotesize{\texttt{arturoscar.lopes@gmail.com}}
\and
Philippe Thieullen\thanks{Partially supported by ANR BLANC07-3\_187245, Hamilton-Jacobi and Weak KAM Theory.}\\
\footnotesize{Institut de Math\'ematiques}\\
\footnotesize{Universit\'e Bordeaux 1, CNRS, UMR 5251}\\
\footnotesize{F-33405 Talence, France}\\
\footnotesize{\texttt{Philippe.Thieullen@math.u-bordeaux1.fr}}
}
\date{\today}

\maketitle

\begin{abstract}
We consider a one-sided transitive subshift of finite type $ \sigma: \Sigma \to \Sigma $ and a H\"older observable $ A $. In the ergodic optimization model,
one is interested in properties of $A$-minimizing probability measures. If $\bar A$ denotes the minimizing ergodic value of $A$, a sub-action $u$ for $A$ is by definition a continuous function such that $A\geq u\circ \sigma-u + \bar A$. We call contact locus of $u$ with respect to $A$ the subset of $\Sigma$ where $A=u\circ\sigma-u + \bar A$. A calibrated sub-action $u$ gives the possibility to construct, for any point $x\in\Sigma$, backward orbits in the contact locus of
$u$. In the opposite direction, a separating sub-action gives the smallest contact locus of $A$, that we call $\Omega(A)$, the set of non-wandering points with respect to $A$.

We prove that separating sub-actions are generic among H\"older sub-actions.
We also prove that, under certain conditions on $\Omega(A)$, any calibrated sub-action is of the form $u(x)=u(x_i)+h_A(x_i,x)$ for some $x_i\in\Omega(A)$, where $h_A(x,y)$ denotes the Peierls barrier of $A$.  We present the proofs in the holonomic optimization model, a formalism which allows to take into account a two-sided transitive subshift of finite type $(\hat \Sigma, \hat \sigma)$.

{\bf To appear in Bull of the Braz. Math. Soc. Vol 40 (4) (2009)}
\end{abstract}

\newpage

\begin{section}{Introduction}

In the {\it ergodic optimization model} (see, for instance, \cite{Bousch1, Bousch2, Bremont, CLT, HY, Jenkinson1, LT1}), given a continuous observable
$ A: X \to \mathbb{R}$, one is interested in understanding which $T$-invariant Borel probability measure $\mu$ of a compact metric space $X$ minimizes
the average $\int_X A\,d\mu$. Such measures are called {\it minimizing probability measures}\footnote{Maximizing probabilities also appear in the
literature. Obviously, replacing the observable $ A $ by $ - A $, both vocabularies can be interchanged and the rephrased statements will be immediately
verified. The maximizing terminology seems more convenient to study the connections with the thermodynamic formalism (see, for example,
\cite{CLT, Leplaideur}).}.

Minimizing probability measures admit dual objects: the {\it sub-actions}. A sub-action $ u: X \to \mathbb{R}$ associated to an observable $A$ enables to
replace $A$ by a cohomologous observable whose ergodic minimizing value is actually the absolute minimum. To each sub-action $u$ one associates a compact subset
of $X$ called {\it contact locus} which contains the support of any minimizing probability measure. A sub-action gives therefore important information on $T$-invariant Borel probability measures that minimize the average of $A$. It is a relevant problem to investigate the existence of a particular sub-action
having the smallest contact locus, that is, the smallest ``trapping region'' of all minimizing probability measures.

In section~\ref{simplified_version}, we give a simplified version for the ergodic optimization model of the main results, namely, of the theorems \ref{principal}, \ref{structure} and \ref{discretestructure}. In section~\ref{basic_concepts}, we recall the definition of the {\it holonomic optimization model} and state the main results. We give in section \ref{proof_theorem_principal} the proof of theorem \ref{principal} and in section \ref{proof_theorem_structure} the proof of theorem \ref{discretestructure}. We address the reader to \cite{GL2} for a proof of theorem~\ref{structure}.
We will adopt throughout the text the point of view which consists in interpreting ergodic optimization problems as questions of variational  dynamics (see, for instance, \cite{CLT, GL2, LT1}), similar to Aubry-Mather technics for Lagrangian systems. For an expository introduction to the general theory of ergodic optimization, we refer the reader to the article of O. Jenkinson (see \cite{Jenkinson1}).

We still would like to point out that one of the main conjectures in the theory of ergodic optimization on compact spaces can be roughly formulated in the
following way: \emph{in any hyperbolic dynamics, a generic H\"older (or Lipshitz) observable possesses an unique minimizing probability measure, which is
supported by a periodic orbit}. Concerning this problem, partial answers were already obtained, among them \cite{Bousch2, CLT, HY, LT1, Morris, TAZ}.
Working with a transitive expanding dynamical system, J. Br\'emont (see \cite{Bremont}) has recently shown how such conjecture might follow from a careful
study of the contact loci of typical sub-actions with finitely many connected components. Of course, such result reaffirms the importance of the study of
sub-actions as well as of their respective contact loci.

In the same dynamical context, we are in particular interested in finding \emph{separating sub-actions}, that is
sub-actions whose contact locus is the smallest one. As mentioned above, these sub-actions give more information
on the minimizing measure(s) than does a general sub-action. Our main theorem (namely, theorem~\ref{principal})
states that such sub-actions are actually generic among the set of H\"older sub-actions. An interesting result we
also present here and which is independent of the previous considerations is an analysis related to the following
situation: it is known that, for each irreducible component of the $A$-non-wandering set, one can associate via the
Peierls barrier a calibrated sub-action. We present in theorem~\ref{discretestructure} sufficient conditions (by no
means necessary) that assure that there exists a dominant one among such calibrated sub-actions.
\medskip

\noindent \textbf{Acknowledgement.} We would like to thanks the referee for his careful reading of our manuscript. This improved the exposition of the final version considerably.

\end{section}

\begin{section}{A simplified version of theorems \ref{principal}, \ref{structure} and \ref{discretestructure}} \label{simplified_version}

Let $ (X, T) $ be a transitive expanding dynamical system, that is, a continuous covering several-to-one  map $ T: X \to X $ on a compact metric space $ X $ whose
inverse branches are uniformly contracting by a factor $ 0 < \lambda < 1 $. We denote by $ \mathcal M_T $ the set of $T$-invariant Borel probability
measures. Our objective in this section is to summarize the conclusions of theorems \ref{principal}, \ref{structure} and \ref{discretestructure} in ergodic optimization theory. We first recall basic definitions from \cite{CLT} (see also \cite{Jenkinson1}).

Given a continuous observable $ A: X \to \mathbb R $, we call {\it ergodic minimizing value} the quantity
\[
\bar A := \min_{\mu \in \mathcal M_T} \int A \; d\mu.
\]
We call {\it $A$-minimizing probability} a measure $ \mu \in \mathcal M_T $ which realizes the above minimum.

We say that a continuous function $ u: X \to \mathbb R $ is a sub-action with respect to the observable $ A $ if
the following inequality holds everywhere on $ X $
\begin{equation*}
A \geq u\circ T - u + \bar A.
\end{equation*}
We would like to emphasize that, although the definition of a sub-action can be extended
to other regularities (for instance, to the class of bounded measurable functions),
we will only consider continuous sub-actions in this paper.

\begin{definition}
A sub-action $u:X\to\mathbb{R}$ is said calibrated if
\[
u(x) = \min_{T(y)=x} [ u(y) + A(y) - \bar A ] \;\; \text{ for all } \; x \in  X.
\]
\end{definition}

\begin{definition}
We call contact locus of a sub-action $u$ the set
\[
\mathbb M_A(u) := (A - u \circ T + u )^{-1}(\bar A).
\]
It is just the subset of $ X $ where $ A = u\circ T - u + \bar A $.
\end{definition}

A point $ x \in X $ is said to be {\it non-wandering with respect to $A$} if, for every $ \epsilon > 0 $, there exists an integer
$ k \ge 1 $ and a point $ y \in X $ such that
$$ d(x, y) < \epsilon, \;\;\; d(x, T^k(y)) < \epsilon \; \text{ and } \;
\Big| \sum_{j = 0}^{k - 1} (A - \bar A) \circ T^j (y) \Big| < \epsilon. $$
We denote by $ \Omega(A) $ the {\it set of non-wandering points} with respect to the observable $ A \in C^0(X) $.
When the observable is H\"older, $\Omega(A)$ is a non-empty compact $T$-invariant set containing the support of all minimizing probability measures. Moreover,
\[
\Omega(A) \subset \bigcap \Big\{ \mathbb M_A(u) \,\,\big| \textrm{ $u$ is a continuous sub-action} \Big\}.
\]
We are interested in finding $u$ so that $\Omega(A)=\mathbb{M}_A(u)$.

\begin{definition}
A sub-action $ u \in C^0(X) $ is said to be separating (with respect to $A$) if it satisfies $ \mathbb M_A(u) = \Omega(A) $.
\end{definition}

The main conclusion of theorem \ref{principal} can be stated in the following way.
The proof of this particular case will not be given and can be adapted from the one of the general situation (see section~\ref{proof_theorem_principal}).

\begin{theorem}\label{principalbis}
Let $ (X, T) $ be a transitive expanding dynamical system on a compact metric space and $ A : X \to \mathbb R $ be a $\theta$-H\"older observable.
Then there exist a $\theta$-H\"older separating sub-action for $ A $. Furthermore, in the $\theta$-H\"older topology, the subset of
$\theta$-H\"older separating sub-actions is generic among all $\theta$-H\"older sub-actions.
\end{theorem}

We will present in theorem~\ref{discretestructurebis} a result of different nature and independent
interest. The item which is totally new on this claim will be item 2.

Contrary to a separating sub-action, a calibrated sub-action $u$ possesses a large contact locus in the sense $T(\mathbb{M}_A(u))=X$. Calibrated sub-actions are built using a particular sub-action called the {\it Peierls barrier}. For H\"older observable $A$, the Peierls barrier of $A$, $ h_A: \Omega(A) \times X \to \mathbb{R}$, is a H\"older calibrated sub-action in the second variable defined by
\begin{multline*}
h_A (x,y) := \lim_{\epsilon \to 0} \; \liminf_{k \to +\infty} \; \inf
\Big\{\sum_{j = 0}^{k - 1} (A - \bar A) \circ T^j (z) \,\big|\,
\\  z \in X,\,\, d(z,x) < \epsilon \textrm{ and } d(T^k(z),y) < \epsilon \Big\}.
\end{multline*}

The equivalent theorem to \ref{structure} may be stated in the following form.

\begin{theorem}\label{structurebis}
Let $ (X, T) $ be a transitive expanding dynamical system on a compact metric space and $ A : X \to \mathbb R $ be a H\"older observable.
Then the set of continuous calibrated sub-actions coincides with the set of functions of the form
\[
u(y)=\min_{x\in\Omega(A)} [\phi(x)+h_A(x,y)],
\quad
\forall\ y\in X,
\]
where $\phi:\Omega(A)\to\mathbb{R}$ is any continuous function satisfying
\[
\phi(y)-\phi(x)\leq h_A(x,y),
\quad
\forall\ x,y\in\Omega(A).
\]
Moreover, $u$ extends $\phi$ and is thus uniquely characterized by $\phi$.
\end{theorem}

The condition $x\sim y \Leftrightarrow h_A(x,y)+h_A(y,x)=0$ defines an equivalent relation on $\Omega(A)$. An equivalence class is called an {\it irreducible component}. It is a closed $T$-invariant set\footnote{We prove these statements in the general setting (see definition-proposition~\ref{relacao de equivalencia}
and proposition~\ref{proposicao componentes}).}.

In the case $\Omega(A)$ is reduced to a finite number of disjoint irreducible components, the set of calibrated sub-actions is parametrized by a finite number of conditions. More precisely, if $\Omega(A)=\sqcup_{i=1}^r C_i$ is equal to a disjoint union of irreducible components and $x_i\in C_i$ are chosen, the {\it sub-action constraint set} is by definition
\[
\mathcal{C}_A(x_1,\ldots,x_r) := \{(u_1,\ldots,u_r)\in\mathbb{R}^r \mid u_j-u_i \leq h_A(x_i,x_j),\quad \forall\ i,j \}.
\]
Therefore, the analogous result to theorem~\ref{discretestructure} can be stated as follows.

\begin{theorem}\label{discretestructurebis}
Let $ (X, T) $ be a transitive expanding dynamical system on a compact metric space and $ A : X \to \mathbb R $ be a H\"older observable. Assume that $\Omega(A)=\sqcup_{i=1}^r C_i$ is equal to a disjoint union of irreducible components.
\begin{enumerate}
\item There is a one-to-one correspondence between the sub-action constraint set and the set of calibrated sub-actions,
\[
\left\{
\begin{array}{l}
(u_1,\ldots,u_r)\in\mathcal{C}_A(x_1,\dots,x_r) \\
{\displaystyle u(x)=\min_{1\leq i \leq r}[u_i+h_A(x_i,x)]}
\end{array}
\right.
\Longleftrightarrow
\left\{
\begin{array}{l}
u \textrm{ is a calibrated sub-action} \\
u_i = u(x_i)
\end{array}
\right..
\]
\item Let $i_0\in\{1,\ldots,r\}$ and $u_{i_0}\in\mathbb{R}$ fixed. Define $u_i=u_{i_0}+h_A(x_{i_0},x_i)$ for all $i$, then $(u_1,\ \ldots,u_r)\in\mathcal{C}_A(x_1,\ldots,x_r)$ and the unique calibrated sub-action $u$ satisfying $u(x_i)=u_i$, for all $i$, is of the form
\[
u(x) := \min_{1\leq i \leq r}[u_i+h_A(x_i,x)] = u_{i_0}+h_A(x_{i_0},x).
\]
\end{enumerate}
\end{theorem}

When the optimizing probability is unique, the calibrated sub-action is unique
(up to additive constants) and generally the proofs of important results are
easiest to discuss.

One of the main issues of the thermodynamic formalism at temperature zero
is the analysis, in the case there are several ergodic maximizing probabilities
for $A$, which of these probabilities the Gibbs states $\mu_{\beta A}$
accumulates, when the inverse temperature parameter $ \beta $ goes to infinite.
It is not clear when there is a unique one in the general H\"older case\footnote{Examples
of Lipschitz observables on the full shift $ \{0,1\}^{\mathbb N} $ for which the zero
temperature limit of the associated Gibbs measures does not exist have been recently
announced (see \cite{CH}).}. In the case of a potential $A$ that depends on finitely
many coordinates, this question is addressed in \cite{Bremont0, Leplaideur}.

Let us denote, in our notation, by $ C_1, C_2, \ldots, C_r $ the different
supports of the ergodic components of the set of maximizing probabilities for
$A$. Then, one can ask: is there a unique one, let us say, with support in $C_{i_0}$
that will be attained as the only limit of Gibbs states $\mu_{\beta A}$ when
$ \beta \to \infty $?

This  question is in some way related to the result of item 2 of theorem~\ref{discretestructurebis}.
Indeed, the dual question can be made for the limits $\frac{1}{\beta} \log \phi_\beta$ when $\beta \to \infty $,
where $\phi_\beta$ is the normalized eigenfunction for the Ruelle operator associated to the potential
$\beta A$. It is well-known that any convergent subsequence will determine a calibrated sub-action,
but is it not clear if there is only one possible limit.

Hence, which one among the various calibrated sub-actions would be chosen?  This is an
important question. All functions of the form  $h_A(x_i,\cdot)$, with $ x_i \in C_i$,
are calibrated sub-actions, for any $ i = 1,2,\ldots,r $. Item 2 of theorem~\ref{discretestructurebis}
gives sufficient conditions to say that a certain $ h_A(x_{i_0},\cdot) $ is preferred in some sense.
We believe this fact is related to the important issues described above.

\end{section}

\begin{section}{Basic Concepts and Main Results}
\label{basic_concepts}

For simplicity, we will restrict the exposition of the {\it holonomic optimization model} to the symbolic dynamics case.
Let $ (\Sigma, \sigma) $ be a one-sided transitive subshift of finite type
given by a $ s \times s $ irreducible transition matrix $ \mathbf M $. More precisely
\[
\Sigma := \Big\{ \mathbf{x} \in \{1, \ldots, s\}^{\mathbb N} \,\big|\, \mathbf{M}(x_j, x_{j+1}) = 1 \text{ for all } j \geq 0 \Big\}
\]
\noindent and $ \sigma $ is the left shift acting on $ \Sigma $ by $ \sigma(x_0, x_1, \ldots) = (x_1, x_2, \ldots) $. Fix
$\lambda \in (0, 1)$. We choose a particular  metric on $\Sigma$ defined by $d(\mathbf x, \bar{\mathbf x}) = \lambda^k$, for any $ \mathbf{x}, \bar{\mathbf{x}} \in \Sigma$, $\mathbf{x} = (x_0, x_1, \ldots) $, $\bar{\mathbf{x}} = (\bar x_0, \bar x_1, \ldots)$ and $ k = \min \{j: x_j \ne \bar x_j \} $.

The holonomic model is a generalization of the ergodic optimization framework. The holonomic model has been introduced first by R. Ma\~n\'e in an attempt to clarify Aubry-Mather theory for continuous time Lagrangian dynamics (see  \cite{CI,Mane}). In this model, the set of invariant minimizing probability measures is replaced by a broader class of measures called holonomic measures. In Aubry-Mather theory for discrete time Lagrangian dynamics on the $n$ dimensional torus  $\mathbb{T}^n$
(see \cite{Gomes}), an holonomic probability measure $\mu(dx,dv)$ is a probability measure on $\mathbb{T}^n \times \mathbb{R}^n $ satisfying
\[
\int_{\mathbb T^n \times \mathbb R^n} f(x + v) \, d\mu(x, v) =
\int_{\mathbb T^n \times \mathbb R^n} f(x) \, d\mu(x, v),
\quad \forall\ f \in C^0(\mathbb T^n),
\]
where the sum $ x + v $ is obviously taken modulo $ \mathbb Z^n $.

One may exploit an interesting analogy with Aubry-Mather theory in symbolic dynamics.
Similarly to the previous example of discrete dynamics, $\Sigma$ will play the role of the ``space of positions'' (analogous to $\mathbb{T}^n$ in the holonomic model) and the set of inverse branches or possible pasts $\Sigma^*$ will play the role of the ``space of immediately anterior velocities'' (analogous to $\mathbb{R}^n$). For a complete exposition and motivation of the holonomic optimization model, see \cite{Garibaldi, GL2}.

We call {\it dual subshift of finite type} the space
\[
\Sigma^* := \Big\{ \mathbf{y} \in \{1, \ldots, s\}^{\mathbb N_*} \,\big|\, \mathbf{M}(y_{j + 1},y_j) = 1 \text{ for all } j \geq 1 \Big\}.
\]
We denote by $\mathbf{y}=(\ldots,y_3,y_2,y_1)$ a point of $\Sigma^*$. We call {\it dual shift} the map $ \sigma^*(\ldots,y_3,y_2, y_1) := (\ldots, y_3, y_2) $. The {\it natural extension} of $(\Sigma, \sigma)$ will play the role of the ``phase space'' (analogous to $\mathbb{T}^n\times\mathbb{R}^n$) and will be identified with a subset of $\Sigma^* \times \Sigma$
\begin{multline*}
\hat \Sigma := \Big\{(\mathbf y, \mathbf x) = (\ldots, y_2, y_1 | x_0, x_1, \ldots) \in \Sigma^* \times \Sigma \,\big|\, \\
\mathbf{x} = (x_0, x_1, \ldots),\,\, \mathbf{y} = (\ldots, y_2, y_1) \textrm{ and } \mathbf{M}(y_1, x_0) = 1 \Big\}.
\end{multline*}
Equivalently, one may write $ \hat \Sigma = \bigcup_{\mathbf x \in \Sigma} \; \Sigma_{\mathbf x}^* \times \{\mathbf x\} $, where
$$ \Sigma_{\mathbf x}^* := \big\{\mathbf y = (\ldots, y_2, y_1) \in \Sigma^*  \,\big|\, \mathbf M(y_1, x_0) = 1 \big\} \quad
\forall \; \mathbf x = (x_0, x_1, \ldots) \in \Sigma. $$

The analogue of the ``discrete Euler-Lagrange map'' is obtained by the usual left shift $\hat{\sigma}$ on the natural extension,
$$ \hat \sigma (\ldots, y_2, y_1 | x_0, x_1, \ldots) = (\ldots, y_1, x_0 | x_1, x_2, \ldots). $$
Consider then $\tau^*: \hat \Sigma \to \Sigma^*$ given by
\[
\tau^*(\mathbf y, \mathbf x) := \tau^*_{\mathbf x}(\mathbf y) := (\ldots,y_2,y_1,x_0).
\]
Notice that $ \tau^*_{\mathbf x}(\mathbf y) \in (\sigma^*)^{-1}(\mathbf y) $.
Similarly, inverse branches of $ \mathbf x \in \Sigma $ with respect to $ \sigma $ are constructed using the map $\tau:\hat{\Sigma}\to\Sigma$,
\[
\tau(\mathbf y, \mathbf x) := \tau_{\mathbf{y}}(\mathbf{x})=(y_1,x_0,x_2,\ldots).
\]
Clearly we have $ \hat\sigma(\mathbf y, \mathbf x) := (\tau^*_{\mathbf{x}}(\mathbf{y}), \sigma(\mathbf{x})) $
and $ \hat{\sigma}^{-1}(\mathbf y, \mathbf x) = (\sigma^*(\mathbf y), \tau_{\mathbf{y}}(\mathbf{x})) $.
Note that $ \tau = \pi \circ \hat \sigma^{-1}$, where $ \pi : \hat \Sigma \to \Sigma$ is the canonical projection onto the $\mathbf{x}$-variable.

Let $ \hat{\mathcal{M}}$ be the set of probability measures over the Borel sigma-algebra  of $ \hat \Sigma $. Instead of considering the set of $\hat{\sigma}$-invariant probability measures\footnote{It is well-known that a H\"older observable defined on the two-sided shift is cohomologous to
an observable that depends just on future coordinates. So a minimization over $\hat\sigma$-invariant probabilities may be reduced to
a minimization over $\sigma$-invariant probabilities.}, we introduce the set of {\it holonomic probability measures},
\[
\hat{\mathcal{M}}_{\textrm{hol}} := \Big\{\hat \mu \in \hat{\mathcal{M}} \,\big|\,
\int_{\hat \Sigma} f(\tau_{\mathbf y}(\mathbf x)) \; d\hat\mu(\mathbf y, \mathbf x) =
\int_{\hat \Sigma} f(\mathbf x) \; d\hat\mu(\mathbf y, \mathbf x), \; \; \forall \; f \in C^0(\Sigma) \Big\}.
\]
It seems important to insist that the holonomic condition demands only the continuous function $ f $ to be defined
on the one-sided shift of finite type $ \Sigma $ and not on the natural extension $ \hat \Sigma $ as would be the case for the characterization of
$\hat\sigma$-invariance.
Observe that $\hat \mu \in \hat{\mathcal{M}}_{\textrm{hol}} $ if, and only if, $ \pi_*(\hat{\mu}) = \pi_*(\hat{\sigma}^{-1}_*(\hat{\mu}))$ if, and only if, $\sigma^{-1}_*(\hat\mu)$ projects onto a $\sigma$-invariant Borel probability measure. As in section~\ref{simplified_version}, we denote by $\mathcal{M}_\sigma$
the set of $\sigma$-invariant Borel probability measures. The triple $ (\hat \Sigma, \hat \sigma, \hat{\mathcal{M}}_{\textrm{hol}}) $ is called the holonomic model. Such a formalism includes the ergodic optimization model discussed in section~\ref{simplified_version} as we will see.

Let $ A \in C^\theta(\hat \Sigma) $ be a  H\"older observable. We would like to emphasize that $ A $ is continuous on the natural extension $\hat \Sigma$.
This is one of the crucial points in the holonomic setting: the possibility of formulating a relevant minimization question for functions
defined on the two-sided shift. Then, we call {\it holonomic minimizing value} of $A$
\begin{eqnarray*}
\bar{A} & := &
\min \Big\{ \int_{\hat \Sigma} A(\mathbf y, \mathbf x) \; d\hat\mu(\mathbf y, \mathbf x) \,\big|\, \hat\mu \in \hat{\mathcal{M}}_{\textrm{hol}} \Big\} \\
& = & \min\Big\{ \int_{\hat\Sigma}A\circ\hat\sigma(\mathbf y, \mathbf x) \, d\hat\mu(\mathbf y, \mathbf x) \,|\, \pi_*(\hat\mu)\in\mathcal{M}_\sigma \Big\}.
\end{eqnarray*}
If $A\circ\hat\sigma=B\circ\pi$ depends only on the $\mathbf{x}$-variable, $\bar A=\bar B$ as in the section \ref{simplified_version}.

The set of {\it minimizing (holonomic) probability measures} is denoted
\[
\hat{\mathcal{M}}_{\textrm{hol}}(A)  :=
\Big\{\hat \mu \in \hat{\mathcal{M}}_{\textrm{hol}} \,\,\big|\, \int_{\hat \Sigma} A(\mathbf y, \mathbf x) \; d\hat\mu(\mathbf y, \mathbf x) = \bar A \Big\}.
\]
A continuous function $u:\Sigma\to\mathbb{R}$ is called {\it sub-action} with respect to $A$ if
\[
u(\mathbf{x}) - u(\tau_{\mathbf{y}}(\mathbf{x})) \leq  A(\mathbf y, \mathbf x) - \bar A,
\quad
\forall\ (\mathbf y, \mathbf x) \in \hat{\Sigma},
\]
or equivalently $ A - \bar A \geq u\circ\pi - u\circ\pi\circ\hat\sigma^{-1}$. We call {\it contact locus} of a sub-action $u$ the set
\[
\hat{\mathbb{M}}_A(u) := (A - u\circ\pi + u\circ\pi\circ\hat\sigma^{-1})^{-1}(\bar A)
\]
where the above inequality becomes an equality, that is, a point $ (\mathbf y, \mathbf x) \in \hat \Sigma $ belongs to $ \hat{\mathbb{M}}_A(u) $
if, and only if, $ u(\mathbf{x}) - u(\tau_{\mathbf{y}}(\mathbf{x})) =  A(\mathbf y, \mathbf x) - \bar A $.
If $A\circ\hat\sigma = B\circ\pi$ for some $B:\Sigma\to\mathbb{R}$, notice that $\pi\circ\hat\sigma^{-1}(\hat{\mathbb{M}}_A(u))=\mathbb{M}_B(u)$.

A calibrated sub-action is a particular sub-action which possesses a large contact locus in the sense that $\pi(\hat{\mathbb{M}}_A(u))=\Sigma$.

\begin{definition}
A sub-action $u:\Sigma\to\mathbb{R}$ is said to be calibrated for $A$ if
\[
u(\mathbf x) = \min_{\mathbf y \in \Sigma_{\mathbf x}^*} \big[ u(\tau_{\mathbf y}(\mathbf x)) + A(\mathbf y, \mathbf x) - \bar A \big],
\quad
\forall \, \mathbf x\in\Sigma,
\]
where recall that $\Sigma_{\mathbf x}^* := \{\mathbf y \in \Sigma^* \,|\, (\mathbf y, \mathbf x) \in \hat{\Sigma} \}$.
\end{definition}

If $\hat B := A\circ\hat\sigma$ and $B(\mathbf{x}) := \min\{\hat B(\mathbf y, \mathbf x) \,|\, \mathbf{y} \in \Sigma_{\mathbf x}^*\}$, then $u$ is a calibrated sub-action for $A$ if, and only if, $u$ is a calibrated sub-action for $B$. Indeed,
\begin{align*}
u(\mathbf{x}) &= \min_{\sigma(\bar{\mathbf{x}})=\mathbf{x}}  \quad  \min_{\mathbf y\in\Sigma_{\mathbf x}^*,\, \tau_{\mathbf{y}}(\mathbf{x})=\bar{\mathbf{x}}}
\big[ u(\bar{\mathbf{x}}) + \hat{B}(\sigma^*(\mathbf{y}),\bar{\mathbf{x}}) - \bar A \big] \\
&= \min_{\sigma(\bar{\mathbf{x}})=\mathbf{x}} \big[ u(\bar{\mathbf{x}}) + B(\bar{\mathbf{x}}) -\bar A \big] =
\min_{\sigma(\bar{\mathbf{x}})=\mathbf{x}} \big[ u(\bar{\mathbf{x}}) + B(\bar{\mathbf{x}}) -\bar B \big].
\end{align*}
(The definition of $B$ gives $\bar B\leq \bar A$ and the calibration gives $\bar B\geq \bar A$.)

A classification theorem for calibrated sub-actions is presented in \cite{GL2}. A central concept is the set of
non-wandering points with respect to $ A $ (previously defined in
\cite{CLT, LT1} in the ergodic optimization model). We call {\it path of length $k$}
a sequence $(\mathbf{z}^0,\ldots,\mathbf{z}^k)$ of points of $\hat{\Sigma}$ such that
\[
\mathbf{z}^{i}=(\mathbf{y}^i,\mathbf{x}^i) \; \textrm{ with } \;
\mathbf{x}^{i}=\tau_{\mathbf{y}^{i + 1}}(\mathbf{x}^{i + 1}),\,\, \forall\ i=0, 1,\ldots,k - 1,
\]
that is, a sequence $(\mathbf{z}^0,\ldots,\mathbf{z}^k)$ where $\mathbf{x}^i=\sigma^i(\mathbf{x}^0)$ for all $i=0,1,\ldots,k$, $\mathbf{x}^0=(x_0,x_1,\ldots,x_{k-1},\mathbf{x}^k)$ and
\begin{multline*}
\mathbf{z}^0=\big(\mathbf{y}^0|x_0,\ldots,x_{k-1},\mathbf{x}^k\big),\,\, \mathbf{z}^1=\big(\sigma^*(\mathbf y^1), x_0|x_1,\ldots,x_{k-1},\mathbf{x}^k\big), \, \ldots,\\
\mathbf{z}^{k-1}=\big(\sigma^*(\mathbf y^{k - 1}), x_{k-2}|x_{k-1},\mathbf{x}^k\big),\,\, \mathbf{z}^k=\big(\sigma^*(\mathbf y^k), x_{k-1}|\mathbf{x}^k\big).
\end{multline*}
Note that the point $ \mathbf{y}^{0} $ is free of any restriction except that $\mathbf{M}(y_1^0,x_0)=1$, more precisely,
one just asks that $ \mathbf y^0 \in \Sigma_{\mathbf x^0}^* $ while $\mathbf{y}^j\in\Sigma_{\mathbf{x}^j}^*\cap(\sigma^*)^{-1}(\Sigma_{\mathbf{x}^{j-1}}^*)$
for $ j = 1, \ldots, k $.
Equivalently, one could present a path in the following way
\begin{multline*}
\mathbf{z}^0=\big(\mathbf{y}^0,\tau_{\mathbf y^1} \circ \tau_{\mathbf y^2} \circ \cdots \circ \tau_{\mathbf y^k}(\mathbf{x}^k)\big),\,\, \mathbf{z}^1=\big(\mathbf{y}^1,\tau_{\mathbf y^2} \circ \cdots \circ \tau_{\mathbf y^k}(\mathbf{x}^k)\big), \, \ldots,\\
\mathbf{z}^{k-1}=\big(\mathbf y^{k - 1}, \tau_{\mathbf y^k}(\mathbf{x}^k)\big),\,\, \mathbf{z}^k=\big(\mathbf y^k,\mathbf{x}^k\big).
\end{multline*}

Given $ \epsilon > 0 $ and $ \mathbf x, \bar{\mathbf x} \in \Sigma $, we say that a path of length $k$, $(\mathbf{z}^0,\ldots,\mathbf{z}^k)$, begins within  $\epsilon $ of $\mathbf{x}$ and ends within $ \epsilon $ of $\bar{\mathbf{x}}$ if $d(\mathbf x^0,\mathbf{x}) < \epsilon$
and $ d(\mathbf{x}^{k},\bar{\mathbf{x}}) < \epsilon $. Denote by $\mathcal{P}_k(\mathbf{x},\bar{\mathbf{x}},\epsilon)$ the set of such paths. Denote by $\mathcal{P}_k(\mathbf{x})$ the set of paths of length $k$ beginning exactly  at $\mathbf{x}$. Notice that a path $ (\mathbf z^0, \ldots, \mathbf z^k) $
belongs to $ \mathcal{P}_k(\mathbf{x}) $ if, and only if, $ \pi(\mathbf z^i) = \sigma^i(\mathbf x) $ for all $ i = 0, 1, \ldots, k $.

A point $ \mathbf x \in \Sigma $ will be called {\it non-wandering with respect to $A$} if, for every $ \epsilon > 0 $,
one can find a path $(\mathbf{z}^0,\ldots,\mathbf{z}^k)$ in $\mathcal{P}_k(\mathbf x, \mathbf x, \epsilon) $, with $ k \ge 1 $, such that
\[
\Big| \sum_{i = 1}^{k} (A - \bar A)(\mathbf{z}^i) \Big| < \epsilon.
\]
We will denote by $ \Omega(A) $ the set of non-wandering
points with respect to $ A $. If $A\circ\hat\sigma=B\circ\pi$, notice that $\Omega(A)=\Omega(B)$ as in section \ref{simplified_version}.

The first two authors have proved in \cite{GL2} that $\Omega(A)$ is a non-empty compact $\sigma$-invariant set and satisfies
\[
\Omega(A) \subset \bigcap \Big\{\pi(\hat{\mathbb{M}}_A(u)) \,\,\big| \textrm{ $u$ is a continuous sub-action} \Big\}.
\]

\begin{remark}
The set $\Omega(A)$ is analogous to the  projected Aubry set in the continuous time Lagrangian dynamics. One could have introduced the corresponding Aubry set $ \hat{\Omega}(A) \subset \hat \Sigma $ and proved
$\pi(\hat{\Omega}(A)) = \Omega(A) $. Unfortunately, even for H\"older observable $A$,
the graph property is not any more true: $\pi : \hat{\Omega}(A) \to \Omega(A)$ is no more bijective. A counter-example can be found in \cite{GL2}. It would be interesting to find the right assumptions on $ A \in C^\theta(\hat \Sigma)$ in order to get this property.
\end{remark}

Contrary to a calibrated sub-action, a {\it separating sub-action} is a sub-action with the smallest contact locus. More precisely,

\begin{definition}
A sub-action $ u \in C^0(\Sigma) $ is said to be separating (with respect to $A$) if it verifies $\pi(\hat{\mathbb{M}}_A(u)) = \Omega(A)$.
\end{definition}

Our first result is the following one.

\begin{theorem}\label{principal}
If $ A : \hat \Sigma \to \mathbb R $ is a $\theta$-H\"older observable, then there exists a $\theta$-H\"older separating
sub-action. Moreover, in the $\theta$-H\"older topology, the subset of $\theta$-H\"older separating sub-actions is generic among all $\theta$-H\"older sub-actions.
\end{theorem}

According to the analogy with continuous time Lagrangian dynamics, sub-actions correspond to viscosity sub-solutions of the stationary Hamilton-Jacobi equation, calibrated sub-actions correspond to the weak KAM solutions introduced by A. Fathi (see \cite{Fathi}) and separating sub-actions correspond to special sub-solutions as described in \cite{FS}.

By adapting the proof of theorem 10 in \cite{GL2} and by using definition \ref{ManePeierlsBarrier} of the the Peierls barrier $h_A$, we obtain a structure theorem for calibrated sub-actions. Such characterization corresponds to the one obtained for weak KAM solutions in Lagrangian dynamics (see \cite{Contreras}).
The proof of the following theorem will be omitted.

\begin{theorem}\label{structure}
Let $A$ be a $\theta$-H\"older observable.
\begin{enumerate}
\item If $u$ is a continuous calibrated sub-action for $A$, then
\[
u(\mathbf x) = \min_{ \bar{\mathbf x} \in \Omega(A)}
\big[ u(\bar{\mathbf x}) + h_A ( \bar{\mathbf x}, \mathbf x) \big].
\]
\item  Conversely, for every continuous application $\phi:\Omega(A)\to\mathbb{R}$ satisfying
\[
\phi(\mathbf{x}) - \phi(\bar{\mathbf{x}}) \leq h_A(\bar{\mathbf{x}},\mathbf{x}),
\quad
\forall\ \mathbf{x},\bar{\mathbf{x}}\in\Omega(A),
\]
\noindent the function $u(\mathbf{x}) := \min_{\bar{\mathbf{x}}\in\Omega(A)} [ \phi(\bar{\mathbf{x}}) + h_A(\bar{\mathbf{x}},\mathbf{x}) ] $
is a continuous calibrated sub-action extending $\phi$ on $\Omega(A)$.
\end{enumerate}
\end{theorem}

In particular, this representation formula for calibrated sub-actions implies immediately that, in order to
compare two such functions, we just need to compare their restrictions to $ \Omega(A) $. For instance, if two
calibrated sub-actions coincide for every non-wandering point with respect to $ A $, then they are the same.

In the case the set of non-wandering points for $A$ is reduced to a finite union of irreducible components $\Omega(A)=C_1\cup\ldots\cup C_r$, the set of calibrated sub-actions admits a simpler characterization. We first show that the condition $\mathbf{x}\sim\bar{\mathbf{x}} \; \Leftrightarrow \;
h_A(\mathbf{x},\bar{\mathbf{x}}) + h_A(\bar{\mathbf{x}},\mathbf{x}) = 0 $ defines an equivalent relation. Each one of its equivalent classes is called
an irreducible component. Let $\bar{\mathbf{x}}^1 \in C_1, \ldots, \bar{\mathbf{x}}^r \in C_r$ fixed. We call {\it sub-action constraint set} the set
\[
\mathcal{C}_A(\bar{\mathbf{x}}^1,\ldots,\bar{\mathbf{x}}^r) = \{(u_1,\ldots,u_r)\in\mathbb{R}^r \mid u_j-u_i \leq h_A(\bar{\mathbf{x}}^i,\bar{\mathbf{x}}^j),\quad \forall\ i,j \}.
\]
\noindent Our second result is the following one.

\begin{theorem}\label{discretestructure}
Let $A$ be a H\"older observable. Assume $\Omega(A)$ is a finite union of disjoint irreducible components, namely, $\Omega(A)=\sqcup_{i=1}^r C_i$.
Let $\bar{\mathbf{x}}^1 \in C_1, \dots, \bar{\mathbf{x}}^r \in C_r$ fixed.
\begin{enumerate}
\item If $u$ is a continuous calibrated sub-action and $u_i := u(\bar{\mathbf{x}}^i)$ for every $i=1,\ldots,r$, then
\[
(u_1,\ldots,u_r) \in \mathcal{C}(\bar{\mathbf{x}}^1,\ldots,\bar{\mathbf{x}}^r)
\quad\textrm{and}\quad
u(\mathbf{x}) = \min_{1\leq i \leq r} \big[ u(\bar{\mathbf{x}}^i) + h_A(\bar{\mathbf{x}}^i,\mathbf{x}) \big].
\]
\item If $(u_1,\ldots,u_r) \in \mathcal{C}(\bar{\mathbf{x}}^1,\ldots,\bar{\mathbf{x}}^r)$ and
$ u(\mathbf{x}) := \min_{1\leq i \leq r} \big[ u_i + h_A(\bar{\mathbf{x}}^i,\mathbf{x}) \big]$,
then $u$ is a continuous calibrated sub-action satisfying $u(\bar{\mathbf{x}}^i)=u_i$ for all $i=1,\ldots,r$.
\item Take $i_0\in\{1,\ldots,r\}$ and $(u_1,\ldots,u_r)$ such that $ u_i := u_{i_0} + h_A(\bar{\mathbf{x}}^{i_0},\bar{\mathbf{x}}^i) $
for all $i=1,\ldots,r$. Then $i_0$ is unique, $(u_1,\ldots,u_r) \in \mathcal{C}(\bar{\mathbf{x}}^1,\ldots,\bar{\mathbf{x}}^r)$ and the unique calibrated sub-action $u$ satisfying $u(\bar{\mathbf{x}}^i)=u_i$, for all $i=1,\ldots,r$, is of the form
$u(\mathbf{x}) = u_{i_0} + h_A(\bar{\mathbf{x}}^{i_0},\mathbf{x})$.
\end{enumerate}
\end{theorem}

The application we present here has a certain similarity to lemma
6 in \cite{AIPS}. We point out that the local character of viscosity solutions (as in definition 1 of \cite{AIPS}) is not present in our setting.

\begin{application}\label{application}
Let $A$ be a H\"older observable. Consider any continuous sub-action $ v $ and a continuous calibrated sub-action $ u $.
\begin{enumerate}
\item Then $u-v$ is constant on every irreducible component and
\[
\min_\Sigma(u-v)=\min_{\Omega(A)}(u-v).
\]
\item Assume $\Omega(A)=\sqcup_{i=1}^r C_i$ is a finite union of disjoint irreducible components.
If $\min_\Sigma(u-v)$ is realized on an unique component $C_{i_1}$ and the other components $C_i$, $i\not=i_1$, are not local minimum for $u-v$, then
\[
u(\mathbf{x})=u(\bar{\mathbf{x}}_{i_1})+h_A(\bar{\mathbf{x}}_{i_1},\mathbf{x}),
\quad
\forall\ \mathbf{x}\in\Sigma,
\]
where $\bar{\mathbf{x}}_{i_1}$ is any point in $C_{i_1}$.
\end{enumerate}
\end{application}

\end{section}

\begin{section}{Proof of theorem~\ref{principal}}
\label{proof_theorem_principal}

We first recall two notions of action potential between two points: the Ma\~n\'e potential and the Peierls barrier. Given $\epsilon > 0$, $ \mathbf x, \bar{\mathbf x} \in \Sigma $ and $k\geq 1$, we denote
\[
S_A^\epsilon(\mathbf x, \bar{\mathbf x}, k) = \inf
\Big\{ \sum_{i = 1}^{k} (A - \bar A)(\mathbf{z}^i) \,\big|\, (\mathbf z^0,\ldots, \mathbf z^{k}) \in \mathcal P_k(\mathbf x, \bar{\mathbf x}, \epsilon) \Big\}.
\]
If $\hat B := A\circ\hat{\sigma}$  and $B := \min\{ \hat B(\mathbf y, \mathbf x) \,|\, \mathbf{y}\in\Sigma_{\mathbf{x}}^* \}$, notice that
\[
S_A^\epsilon(\mathbf x, \bar{\mathbf x}, k) = \inf
\Big\{ \sum_{i=0}^{k-1}(B-\bar B)\circ\sigma^i(\mathbf{x}^0) \,\big|\, d(\mathbf{x}^0,\mathbf{x})<\epsilon, \,\, d(\sigma^k(\mathbf{x}^0),\bar{\mathbf{x}})<\epsilon \Big\}.
\]

\begin{definition}\label{ManePeierlsBarrier}
\noindent We call Ma\~n\'e potential the function $ \phi_A: \Sigma \times \Sigma \to \mathbb R \cup \{+ \infty \} $ defined by
\[
\phi_A(\mathbf x, \bar{\mathbf x}) = \lim_{\epsilon \to 0} \, \inf_{k \geq 1} \, S_A^\epsilon(\mathbf x, \bar{\mathbf x}, k).
\]
\noindent We call Peierls barrier the function $ h_A: \Sigma \times \Sigma \to \mathbb R \cup \{+ \infty \} $ defined by
\[
h_A(\mathbf x, \bar{\mathbf x}) = \lim_{\epsilon \to 0} \, \liminf_{k \to +\infty} \,  S_A^\epsilon(\mathbf x, \bar{\mathbf x}, k).
\]
\end{definition}

Clearly, $ \phi_A \le h_A $ and both functions are lower semi-continuous. We summarize the main properties of these action potentials.

\begin{proposition}\label{propriedadesbasicas}
Let $A$ be a H\"older observable.
\begin{enumerate}
\item If $u$ is  a continuous sub-action, then $u(\bar{\mathbf x})- u(\mathbf x) \leq  \phi_A(\mathbf x, \bar{\mathbf x})$.
\item For any points $ \mathbf x, \bar{\mathbf x}, \bar{\bar{\mathbf x}} \in \Sigma $,
$ \phi_A(\mathbf x, \bar{\bar{\mathbf x}}) \le \phi_A(\mathbf x, \bar{\mathbf x}) + \phi_A(\bar{\mathbf x}, \bar{\bar{\mathbf x}}) $.
\item Given a point $ \mathbf x \in \Sigma $, if there exists a positive integer $ L $ such that
$ 0 < L < \min \{j > 0 : \sigma^j(\mathbf x) = \mathbf x\} \le + \infty $, then
\[
\phi_A(\mathbf x, \mathbf x) = \phi_A(\mathbf x, \sigma^L(\mathbf x)) + \phi_A(\sigma^L(\mathbf x), \mathbf x).
\]
Moreover, if $\phi_A(\mathbf x, \mathbf{x})<+\infty$, then there exists a path of length $L$,
$(\bar{\mathbf{z}}^0 = (\bar{\mathbf y}^0, \bar{\mathbf x}^0), \ldots, \bar{\mathbf{z}}^L = (\bar{\mathbf y}^L, \bar{\mathbf x}^L))$,
beginning at $\mathbf{x}$ $(\bar{\mathbf{x}}^j=\sigma^j(\mathbf{x})$ for all $j = 0,\ldots, L)$, such that
\[
\phi_A(\mathbf x, \sigma^L(\mathbf x)) = \sum_{j = 1}^L (A - \bar A) (\bar{\mathbf{z}}^j).
\]
\item For any points $ \mathbf x, \bar{\mathbf x}, \bar{\bar{\mathbf x}} \in \Sigma $ and any sequence $ \{\bar{\mathbf  x}^l \} $ converging to $ \bar{\mathbf x} $,
\[
h_A(\mathbf x, \bar{\bar{\mathbf x}}) \leq \liminf_{l \to +\infty} \phi_A(\mathbf x, \bar{\mathbf x}^l) + h_A(\bar{\mathbf x}, \bar{\bar{\mathbf x}}).
\]
\item If $\mathbf{x}\in\Sigma$, then
$ \mathbf x \in \Omega(A) \, \Leftrightarrow \, \phi_A(\mathbf x, \mathbf x) = 0 \, \Leftrightarrow \, h_A(\mathbf x, \mathbf x) = 0 $.
\item If $ \mathbf x \in \Omega(A) $, then $ \phi_A(\mathbf x, \cdot) = h_A(\mathbf x, \cdot) $ and $h_A(\mathbf x, \cdot)$ is a H\"older calibrated sub-action with respect to the second variable.
\end{enumerate}
\end{proposition}

This proposition shows how to construct H\"older calibrated sub-actions without the use of the Lax-Oleinik fixed point method.

\begin{remark}
In Lagrangian Aubry-Mather theory on a compact manifold $ M $, it is well known that, for any point  $ x \in M $, the map $ y \in M \mapsto h(x, y) \in \mathbb R $ defines a weak KAM solution, where $ h : M \times M \to \mathbb R $ denotes the  corresponding Peierls barrier. The analogous result for $ h_A(\mathbf x, \cdot) $ is however false in the holonomic optimization model. Using item 3, it is not difficulty to built examples where
\[
\lim_{L \to + \infty} \phi_A(\mathbf x, \sigma^L(\mathbf x)) = \lim_{L \to + \infty} h_A(\mathbf x, \sigma^L(\mathbf x)) = + \infty,
\]
which shows that $ h_A(\mathbf x, \cdot) $ is not always a continuous function.
\end{remark}

\begin{proof}[Proof of proposition~\ref{propriedadesbasicas}]
Items 1, 2, 5 and 6 are well known and a demonstration can be found, for instance, in \cite{CLT, GL2}.
So let us prove items 3 and 4.

\medskip

\noindent {\it Item 3.} We already know from item 2, that
\[
\phi_A(\mathbf{x},\mathbf{x}) \leq \phi_A(\mathbf{x},\sigma^L(\mathbf{x})) + \phi_A(\sigma^L(\mathbf{x}),\mathbf{x}).
\]
\noindent  Define $ \eta = \min \{d(\sigma^i(\mathbf x), \sigma^j(\mathbf x)) : 0 \le i < j \le L \} $. Fix $ \gamma > 0 $ and take
$ \epsilon \in (0, \min\{\lambda, \eta/2\}) $ such that $\textrm{H\"old}(A) L \epsilon^\theta < \gamma$.
Consider also $ \rho \in (0, \epsilon) $ such that $ d(\mathbf x, \bar{\mathbf x}) < \rho $ implies
$ d(\sigma^j(\mathbf x), \sigma^j(\bar{\mathbf x})) < \epsilon $ for $ 1 \le j \le L $.  Take then a path $ (\mathbf{z}^0,\ldots,\mathbf{z}^l) \in \mathcal P_l(\mathbf x, \mathbf x, \rho) $ satisfying
\[
\sum_{j = 1}^{l} (A - \bar A)(\mathbf{z}^j) < \inf_{k \geq 1} S_A^\rho(\mathbf x, \mathbf x, k) + \gamma \leq \phi_A(\mathbf{x},\mathbf{x})+\gamma.
\]
Let $\mathbf{z}^j=(\mathbf{y}^j, \mathbf{x}^j)$ where $\mathbf{x}^j=\sigma^j(\mathbf{x}^0)$ for all $j=0,1,\ldots,l$.
We claim that $ l > L $. Indeed, $\rho$ has been chosen so that, for each $j \in \{1, 2, \ldots, L\}$,
\[
d(\mathbf x^{j}, \mathbf x) = d(\sigma^j(\mathbf x^0), \mathbf x) \ge  d(\sigma^j(\mathbf x), \mathbf x) - d(\sigma^j(\mathbf x), \sigma^j(\mathbf x^0))
> \eta - \epsilon > \rho.
\]

Introduce a new path $(\bar{\mathbf{z}}^0,\ldots,\bar{\mathbf{z}}^L)\in\mathcal{P}_L(\mathbf{x},\mathbf{x},\epsilon)$ given by
$ \bar{\mathbf{z}}^j=(\mathbf{y}^j,\sigma^j(\mathbf{x})) $, for all $ j=0,\ldots, L $.
The definition of $\rho$ guarantees
\[
\sum_{j = 1}^{L} (A - \bar A)(\bar{\mathbf{z}}^j)
< \sum_{j = 1}^{L} (A - \bar A)(\mathbf{z}^j) +  \text{H\"old}_\theta(A) L \epsilon^\theta \leq \sum_{j = 1}^{L} (A - \bar A)(\mathbf{z}^j) + \gamma.
\]
Notice that $ (\mathbf{z}^L,\ldots,\mathbf{z}^l)\in\mathcal{P}_{l-L}(\sigma^L(\mathbf{x}),\mathbf{x},\epsilon) $. We finally obtain
\begin{align*}
\inf_{k\geq 1} S_A^\epsilon(\mathbf x,\sigma^L(\mathbf x), k) &+ \inf_{k\geq 1} S_A^\epsilon(\sigma^L(\mathbf x), \mathbf x, k) \\
&\leq \sum_{j = 1}^{L} (A - \bar A)(\bar{\mathbf{z}}^j) + \sum_{j = L+1}^{l} (A - \bar A)(\mathbf{z}^j) \\
&\leq \sum_{j = 1}^{L} (A - \bar A)(\mathbf{z}^j) + \sum_{j = L+1}^{l} (A - \bar A)(\mathbf{z}^j) + \gamma \\
&\leq \inf_{k \geq 1} S_A^\rho(\mathbf x, \mathbf x, k) + 2\gamma \leq \phi_A(\mathbf{x},\mathbf{x})+2\gamma.
\end{align*}
By letting $\epsilon$ goes to $0$ and $\gamma$ to $0$, we get
\[
\phi_A(\mathbf x, \sigma^L(\mathbf x)) + \phi_A(\sigma^L(\mathbf x), \mathbf x) \leq \phi_A(\mathbf x, \mathbf x).
\]

The first part of item 3 is proved. To prove the second part, the previous computation shows that, for any sufficiently small $\epsilon$, there exists a path $(\bar{\mathbf{z}}^0_\epsilon,\ldots,\bar{\mathbf{z}}^L_\epsilon)\in\mathcal{P}_L(\mathbf{x})$ such that
\[
\sum_{j = 1}^{L} (A - \bar A)(\bar{\mathbf{z}}^j_\epsilon) + \inf_{k\geq 1} S_A^\epsilon(\sigma^L(\mathbf x), \mathbf x, k) \leq \phi_A(\mathbf{x},\mathbf{x})+2\gamma.
\]
By taking accumulation points of $\bar{\mathbf{z}}^j_\epsilon$ when $\epsilon\to0$, we obtain, for any $\gamma$, a path $(\bar{\mathbf{z}}^0,\ldots,\bar{\mathbf{z}}^L)$ such that
\[
\phi_A(\mathbf x, \sigma^L(\mathbf x)) \leq \sum_{j = 1}^{L} (A - \bar A)(\bar{\mathbf{z}}^j) \leq \phi_A(\mathbf{x},\mathbf{x}) - \phi_A(\sigma^L(\mathbf x), \mathbf x) + 2\gamma
\]
\noindent The result follows from item 2 and by taking once more accumulation points of $\bar{\mathbf{z}}^j$ when $\gamma\to0$.

\medskip

\noindent {\it Item 4.} Since $ \phi_A $ is lower semi-continuous, the statement is equivalent to
\[
h_A(\mathbf x, \bar{\bar{\mathbf x}}) \leq  \phi_A(\mathbf x, \bar{\mathbf x}) + h_A(\bar{\mathbf x}, \bar{\bar{\mathbf x}}),
\quad
\forall\  \mathbf x, \bar{\mathbf x}, \bar{\bar{\mathbf x}} \in \Sigma.
\]

Fix $\gamma>0$ and $\epsilon \in (0, \lambda/2) $ such that $\textrm{H\"old}(A)(2\epsilon)^\theta/(1-\lambda^\theta)<\gamma$. There exists a path  $(\mathbf{z}^0,\ldots,\mathbf{z}^k) \in \mathcal{P}_k(\mathbf x, \bar{\mathbf x}, \epsilon)$ such that
\[
\sum_{j = 1}^{k} (A - \bar A)(\mathbf z^j) < \inf_{n \geq 1} S_A^\epsilon(\mathbf x, \bar{\mathbf x}, n) + \gamma.
\]
For any $ N \geq 1 $, there exists a path $(\bar{\mathbf z}^0,\ldots,\bar{\mathbf z}^{l})
\in \mathcal P_{l}(\bar{\mathbf x}, \bar{\bar{\mathbf x}}, \epsilon) $ of length $l\geq N$ such that
\[
\sum_{j = 1}^{l} (A - \bar A)(\bar{\mathbf z}^j) < \inf_{n \ge N} S_A^\epsilon(\bar{\mathbf x}, \bar{\bar{\mathbf x}}, n) + \gamma.
\]

We define a path $(\bar{\bar{\mathbf{z}}}^0,\ldots,\bar{\bar{\mathbf{z}}}^{k+l})\in\mathcal{P}_{k+l}(\mathbf x, \bar{\bar{\mathbf x}}, 3\epsilon)$ in the following way
\begin{gather*}
\bar{\bar{\mathbf{z}}}^j=\bar{\mathbf{z}}^{j-k},\,\,
\forall\ j=k+1,\ldots,k+l,
\quad\quad
\bar{\bar{\mathbf{z}}}^j=(\bar{\bar{\mathbf{y}}}^j,\bar{\bar{\mathbf{x}}}^j),\,\,
\forall\ j=0,\ldots,k, \\
\bar{\bar{\mathbf{y}}}^j=\mathbf{y}^j,\,\,
\forall\ j=0,\ldots,k,
\quad\quad
\bar{\bar{\mathbf{x}}}^k=\bar{\mathbf{x}}^0,\,\,
\bar{\bar{\mathbf{x}}}^{j-1}=\tau_{\mathbf{y}^j}(\bar{\bar{\mathbf{x}}}^j),\,\,
\forall\ j=1,\ldots,k.
\end{gather*}
\noindent We notice that $d(\bar{\bar{\mathbf{x}}}^j,\mathbf{x}^j) \leq \lambda^{k-j}d(\bar{\bar{\mathbf{x}}}^k,\mathbf{x}^k)$, for all $j=0,\ldots,k$. Since
\[
d(\bar{\bar{\mathbf{x}}}^k,\mathbf{x}^k)= d(\bar{\mathbf{x}}^0,\mathbf{x}^k) \leq d(\bar{\mathbf{x}}^0,\bar{\mathbf{x}}) +d(\bar{\mathbf{x}},\mathbf{x}^k)< 2\epsilon,
\]
\noindent we obtain $d(\bar{\bar{\mathbf{x}}}^0,\mathbf{x}) \leq \lambda^k 2\epsilon +\epsilon < 3\epsilon$. Hence, it follows that
\begin{align*}
\inf_{n \ge N} S_A^{3\epsilon}(\mathbf x, \bar{\bar{\mathbf x}}, n)
&\leq \sum_{j = 1}^{k+l} (A - \bar A)(\bar{\bar{\mathbf z}}^j) \\
&\leq \sum_{j = 1}^{l} (A - \bar A)(\bar{\mathbf z}^j) + \sum_{j = 1}^{k} (A - \bar A)(\mathbf z^j) + \frac{(2\epsilon)^\theta}{1 - \lambda^\theta}\text{H\"old}_\theta(A) \\
&\leq \inf_{n \geq 1} S_A^\epsilon(\mathbf x, \bar{\mathbf x}, n) + \inf_{n \ge N} S_A^\epsilon(\bar{\mathbf x}, \bar{\bar{\mathbf x}}, n) + 3\gamma \\
&\leq \phi_A(\mathbf x, \bar{\mathbf x}) + \inf_{n \ge N} S_A^\epsilon(\bar{\mathbf x}, \bar{\bar{\mathbf x}}, n) + 3\gamma.
\end{align*}
\noindent By taking first $N\to+\infty$, then $\epsilon\to 0$ and $\gamma\to 0$, we get
\[
h_A(\mathbf x, \bar{\bar{\mathbf x}}) \leq \phi_A(\mathbf x, \bar{\mathbf x}) + h_A(\bar{\mathbf x}, \bar{\bar{\mathbf x}}).
\]
\end{proof}

Other properties of the Ma\~n\'e potential and the Peierls barrier can be derived from the previous proposition. For instance, item 4 gives us the following inequality
\[
h_A(\mathbf x, \bar{\bar{\mathbf x}}) \le h_A(\mathbf x, \bar{\mathbf x}) + h_A(\bar{\mathbf x}, \bar{\bar{\mathbf x}}),
\quad
\forall\ \mathbf x, \bar{\mathbf x}, \bar{\bar{\mathbf x}} \in \Sigma.
\]

We now begin the proof of theorem \ref{principal}. It follows immediately from the next lemma.

\begin{lemma}\label{secundario}
Let $ D \subset \Sigma $ be an open set containing $ \Omega(A) $. Denote by
$\mathcal{D}_A$ the subset of H\"older sub-actions $ u $ such that $ \pi(\mathbb M_A(u)) \subset D $. Then, for the H\"older topology, $\mathcal{D}_A$ is an open dense subset of the H\"older sub-actions.
\end{lemma}

We only need a few lines to show that lemma~\ref{secundario} yields theorem~\ref{principal}.
As a matter of fact, if one considers, for each positive integer $ j $, the open set
$ D_j = \{ \mathbf x \in \Sigma \,|\, d(\mathbf x, \Omega(A)) < 1/j \} $ and the corresponding open dense subset of H\"older sub-actions $ \mathcal{D}_{A,j} $, then the set of H\"older separating sub-actions contains the countable intersection $\cap_{j>0}\mathcal{D}_{A,j}$.

\begin{proof}[Proof of lemma \ref{secundario}] We only discuss the denseness of $ \mathcal{D}_A$.

\medskip
\noindent {\it Part 1.} Let $ v $ be any H\"older sub-action for $ A $. We will show that, for every $ \mathbf x \notin D $, there exists a H\"older sub-action
$ v_{\mathbf x} $ as close as we want to $ v $ in the H\"older topology with a projected  contact locus disjoint from $ \mathbf x $, that is, $ \mathbf x \notin \pi(\mathbb M_A(v_{\mathbf x})) $ or
\[
v_{\mathbf x} (\mathbf x) - v_{\mathbf x}(\tau_{\mathbf y}(\mathbf x)) < A(\mathbf y, \mathbf x) - \bar A,
\quad
\forall\ \mathbf y \in \Sigma_{\mathbf x}^*.
\]

Let $ \mathbf x \notin D $. We discuss two cases.

\medskip
\noindent {\it Case a.} We assume there exists an integer $k \geq 0$ such that, for every path of length $k$ beginning at $\mathbf{x}$,
$(\mathbf{z}^0 = (\mathbf y^0, \mathbf x), \ldots, \mathbf{z}^k =(\mathbf y^k, \sigma^k(\mathbf x))) \in \mathcal{P}_k(\mathbf{x})$,
the terminal point $\mathbf{z}^k \not\in \mathbb{M}_A(v)$. If $k=0$, we choose $v_{\mathbf{x}}=v$. Assume now $k\geq 1$. Let
\[
B := A - \bar A - v \circ \pi + v \circ \pi \circ \hat \sigma^{-1} \geq 0
\]
be the associated normalized observable ($B\geq0$ and $\bar B=0$).
We recall that $ \tau_{\mathbf{y}^j}(\sigma^{j}(\mathbf x)) = \sigma^{j-1}(\mathbf x) $, for all $j=1,\ldots,k$.
So by hypothesis
\begin{equation}
B(\mathbf{z}^k)=B(\mathbf y^k, \sigma^k(\mathbf{x})) > 0,
\;\; \forall\ \mathbf{y}^k \in \Sigma_{\sigma^{k}(\mathbf x)}^*
\textrm{ s.t. } \sigma^{k-1}(\mathbf x)=\tau_{\mathbf{y}^k }(\sigma^{k}(\mathbf x)). \tag{I}
\end{equation}

Notice first that, if $(\bar{\mathbf{z}}^0,\ldots,\bar{\mathbf{z}}^k)$ is a path of length $k$ and $ \gamma \in (0,1) $ is any constant, as $ B $ is non-negative,
one has
\begin{align}
B(\bar{\mathbf{z}}^0) &= \sum_{j = 0}^{k-1} B(\bar{\mathbf z}^j) - \sum_{j = 1}^{k} B(\bar{\mathbf z}^j) + B(\bar{\mathbf z}^k) \nonumber \\
&\geq \gamma \sum_{j = 0}^{k-1}  B(\bar{\mathbf z}^j) - \gamma  \sum_{j = 1}^{k} B(\bar{\mathbf z}^j) + \gamma B(\bar{\mathbf z}^k)
\geq \gamma \sum_{j = 0}^{k-1} B(\bar{\mathbf z}^j) - \gamma \sum_{j = 1}^{k}  B(\bar{\mathbf z}^j). \tag{II}
\end{align}

Let $w_{k} : \Sigma \to \mathbb R$ be the function given by
\[
w_k(\bar{\mathbf x}) :=
\inf \Big\{ \sum_{j = 1}^{k} B(\bar{\mathbf z}^j) \,\big|\, (\bar{\mathbf z}^0, \ldots,\bar{\mathbf z}^{k}) \in \mathcal P_k(\bar{\mathbf{x}}) \Big\},
\quad
\forall\ \bar{\mathbf x} \in \Sigma.
\]
Because $ \mathcal P_k(\bar{\mathbf{x}}) $ is a closed subspace of the compact space $ \hat \Sigma^{k + 1} $, the above infimum is effectively a minimum.
Moreover, since the application
$ C(\bar{\mathbf{x}}) := \min \{ B \circ \hat \sigma(\bar{\mathbf y},\bar{\mathbf{x}}) \,|\, \bar{\mathbf{y}}\in\Sigma_{\bar{\mathbf{x}}}^* \}$ is H\"older,
$w_k = \sum_{j=0}^{k-1} C\circ\sigma^j$ is also H\"older\footnote{We leave the details to the reader. In particular, one shall note that
$w_k = \sum_j C\circ\sigma^j$ means $ \min(a + b) = \min a + \min b $, which indicates the importance of what a path is.}.

We first prove that $-\gamma w_k$ is a sub-action. Let $\bar{\mathbf{x}}\in\Sigma$ and $\bar{\mathbf{y}}\in\Sigma_{\bar{\mathbf{x}}}^*$. There exists a path of length $k$, $ (\bar{\mathbf z}^0,\ldots,\bar{\mathbf z}^k)$, beginning at $\bar{\mathbf{x}}$ and realizing the minimum
\[
 w_k(\bar{\mathbf x}) = \sum_{j = 1}^{k} B(\bar{\mathbf z}^j).
\]
Notice the only constraint on $\bar{\mathbf{y}}^0$ is $\bar{\mathbf{y}}^0\in\Sigma_{\bar{\mathbf{x}}}^*$, besides $\bar{\mathbf{y}}^0$ does not appear in the previous sum. Choose $\bar{\mathbf{y}}^0=\bar{\mathbf{y}}$, $\bar{\mathbf{y}}^{-1}\in\Sigma_{\bar{\mathbf{x}}^{-1}}^*$ and call $\bar{\mathbf{x}}^{-1}=\tau_{\bar{\mathbf{y}}}(\bar{\mathbf{x}})$. Then $(\bar{\mathbf{z}}^{-1},\bar{\mathbf{z}}^0,\ldots,\bar{\mathbf{z}}^{k-1})$ is a path of length $k$ beginning at $\tau_{\bar{\mathbf{y}}}(\bar{\mathbf{x}})$. So denote $ \bar{\mathbf z} := (\bar{\mathbf x}, \bar{\mathbf y}) $. Thanks to inequality (II)
\[
B(\bar{\mathbf{z}})=B(\bar{\mathbf{z}}^0) \geq \gamma \sum_{j = 0}^{k-1} B(\bar{\mathbf z}^j) - \gamma \sum_{j = 1}^{k}  B(\bar{\mathbf z}^j) \ge
\gamma w_k(\tau_{\bar{\mathbf{y}}}(\bar{\mathbf{x}})) - \gamma w_k(\bar{\mathbf{x}}),
\]
which shows $-\gamma w_k$ is a sub-action for $ B $.
Moreover, given any $\mathbf{y}\in\Sigma_{\mathbf{x}}^*$, the same computation for $ \mathbf{z} := (\mathbf x, \mathbf y) $ instead of
$\bar{\mathbf{z}}$ and (I) assure that
\[
B(\mathbf{z}) - \gamma w_k(\tau_{\mathbf{y}}(\mathbf{x})) + \gamma w_k(\mathbf{x}) \geq \gamma B(\mathbf{z}^k) > 0.
\]
We have proved that $ \mathbf x \notin \pi(\mathbb M_B(-\gamma w_k))=\pi(\mathbb{M}_A(v-\gamma w_k))$.

Since $ \gamma $ can be taken as small as we want, we have shown the existence of a H\"older sub-action $v_{\mathbf x} = v - \gamma w_{k} $ close to $v$ in the H\"older topology satisfying $ \mathbf x \notin \pi(\mathbb M_A(v_{\mathbf x})) $.

\medskip

\noindent {\it Case b.} We suppose that, for every integer $ k\geq0 $, one can find a path of length $k$,
$(\mathbf{z}^0,\ldots,\mathbf{z}^{k})$, beginning at $\mathbf{x}$, such that $\mathbf{z}^k \in \mathbb M_A(v)$,
or equivalently $ B(\mathbf{z}^k) = 0 $ with $B$  as before. In other words, there exists $ \mathbf y^0 \in \Sigma_{\mathbf x}^* $ with
$ B(\mathbf y^0, \mathbf x) = 0 $ and, for any $k\geq1$, there exists $\mathbf{y}^k\in\Sigma_{\mathbf{x}^k}^*\cap(\sigma^*)^{-1}(\Sigma_{\mathbf{x}^{k-1}}^*)$ such that $B(\mathbf{y}^k,\mathbf{x}^k)=0$, where $\mathbf{x}^k=\sigma^k(\mathbf{x})$. Define $\bar{\mathbf{z}}^0=(\mathbf{y}^0,\mathbf{x})$ and
$\bar{\mathbf{z}}^k=(\mathbf{y}^k,\mathbf{x}^k)$ for all $k\geq1$. Notice that $(\bar{\mathbf{z}}^0,\ldots,\bar{\mathbf{z}}^k)$ is now a path of arbitrary length
$ k $, beginning at $\mathbf{x}$, which satisfies $B(\bar{\mathbf{z}}^j)=0$ for $j = 0, \ldots, k $.

Let $\bar{\mathbf{x}}\in\Omega(A)=\Omega(B)$ be any limit point of $(\mathbf{x}^k)_k$ chosen once for all. Let $ w := h_B(\bar{\mathbf x},\cdot)$ be the H\"older sub-action for $ B $ given by the corresponding Peierls barrier. Notice that by the definition of the Peierls barrier (see definition~\ref{ManePeierlsBarrier})
we clearly get $ h_B \ge 0 $, since $B\geq0$ and $\bar B=0$. Furthermore, we remark that $ \phi_B (\mathbf x, \sigma^k(\mathbf x)) = 0 $ for all $k\geq1$ and that
\[
w(\mathbf{x}) = h_B(\bar{\mathbf x}, \mathbf x) = \liminf_{k \to +\infty} \phi_B(\mathbf x, \sigma^k(\mathbf x)) + h_B(\bar{\mathbf x}, \mathbf x)
\ge h_B(\mathbf x, \mathbf x) > 0.
\]
Here we have used item 4 of proposition~\ref{propriedadesbasicas} to obtain the first inequality and item 5 of the same proposition to assure
the strict inequality since $ \mathbf x \notin D \supset \Omega(A) = \Omega(B) $.

Let $ \gamma \in  (0,1) $ be any real number as close to 0 as we want. We claim that $ \mathbf x $ satisfies again the first case, namely, there exists $k\geq1$ such that, for any path of length $k$, $(\mathbf{z}^0 = (\mathbf y^0,\mathbf x^0), \ldots, \mathbf{z}^k = (\mathbf y^k, \mathbf x^k))$, beginning at $\mathbf{x}$,
one has
\[
B(\mathbf{z}^k) - \gamma h_B(\bar{\mathbf{x}},\mathbf{x}^k) +\gamma h_B(\bar{\mathbf{x}},\mathbf{x}^{k-1}) > 0.
\]
(Notice that $\gamma w$ is again a sub-action for $B$ since $B$ is non-negative.) Indeed, by contradiction, for any integer $ k \ge 0 $, we would have
a path of length $ k $, $ (\mathbf{z}^0 = (\mathbf y^0,\mathbf x^0), \ldots, \mathbf{z}^k = (\mathbf y^k, \mathbf x^k)) $,
beginning at $ \mathbf x $, such that $ \mathbf z^k \in \mathbb M_B(\gamma w) $, which would yield
\[
0 \leq B(\mathbf{z}^k) = \gamma h_B(\bar{\mathbf{x}},\mathbf{x}^k) - \gamma h_B(\bar{\mathbf{x}},\mathbf{x}^{k-1}), \quad \forall \, k \ge 1.
\]
On the one hand, from the inequality $ \gamma h_B(\bar{\mathbf{x}},\mathbf{x}^{k-1}) \le \gamma h_B(\bar{\mathbf{x}},\mathbf{x}^k) $,
we would obtain $0 < w(\mathbf x) = h_B(\bar{\mathbf{x}},\mathbf{x}) \leq h_B(\bar{\mathbf{x}},\mathbf{x}^k)$ for all $k\geq1$.
On the other hand, by taking a subsequence of $\{\mathbf{x}^k\} = \{\sigma^k(\mathbf x)\} $ converging to $ \bar{\mathbf x} $, $h_B(\bar{\mathbf{x}},\mathbf{x}^k)$ would converge to $ h_B(\bar{\mathbf{x}},\bar{\mathbf{x}}) = 0 $, since $ \bar{\mathbf x} \in \Omega(B) $. We have thus obtained a contradiction.
Hence, case (a) implies that there exists a sub-action $v_{\mathbf{x}}$, close to $v$ in the H\"older topology, satisfying
$ \mathbf x \notin \pi(\mathbb M_A(v_{\mathbf x})) $.

\medskip

\noindent{\it Part 2.} We have just proved that, for any $ \mathbf x \notin D $, there exists a sub-action $ v_{\mathbf x} $ close to $ v $ and a
ball $ B(\mathbf x, \epsilon_{\mathbf x}) $ of radius $ \epsilon_{\mathbf x} > 0 $ centered at $ \mathbf x $
such that
$$ \forall \; \bar{\mathbf x} \in B(\mathbf x, \epsilon_{\mathbf x}), \quad \bar{\mathbf x} \notin \pi(\mathbb M_A(v_{\mathbf x})). $$
We can extract from the family of these balls $ \{ B(\mathbf x, \epsilon_{\mathbf x}) \}_{\mathbf x} $ a finite family indexed by
$ \{ \mathbf x^j \}_{1 \le j \le K} $ which is still a covering of the compact set $ \Sigma \setminus D $. Let
$$ u = \frac{1}{K} \sum_{j = 1}^K v_{\mathbf x^j}. $$
Then it is easy to check that $ u $ is a H\"older sub-action for $ A $ satisfying $ \pi(\mathbb M_A(u)) \subset D $, namely, $ u \in \mathcal D_A $.
Since each sub-action $ v_{\mathbf x} $ can be taken as close as we want to $ v $ in the H\"older topology, the same is true for $ u $.
\end{proof}

\end{section}

\begin{section}{Proof of theorem \ref{discretestructure}}
\label{proof_theorem_structure}

It was proved in \cite{GL2} that the projection of the support of a minimizing probability measure $\hat\mu$ is included into the $A$-non-wandering set $\Omega(A)$
when such projection is ergodic.
If $\pi_*\hat\mu$ is ergodic, $\pi(\textrm{supp}(\hat\mu))$ may be seen as an irreducible component in the sense that any two points can be joined by an $\epsilon$-closed trajectory. We introduce here a more general notion of irreducibility.

\begin{definition-proposition}\label{relacao de equivalencia}
Let $A : \hat \Sigma \to \mathbb R $ be a H\"older observable. We say that two points $\mathbf{x},\bar{\mathbf{x}}$ of $\Omega(A)$
are equivalent and write $\mathbf{x}\sim\bar{\mathbf{x}}$ if
\[
h_A(\mathbf{x},\bar{\mathbf{x}}) + h_A(\bar{\mathbf{x}},\mathbf{x}) = 0.
\]
Then $\sim$ is an equivalent relation. Its equivalent classes are called irreducible components.
\end{definition-proposition}

\begin{proof}
It is obvious that $\sim$ is reflexive ($h_A(\mathbf{x},\mathbf{x})=0 \Leftrightarrow \mathbf{x}\in\Omega(A)$) and symmetric. Let $u$ be a continuous sub-action and $B:=A-\bar A-u\circ\pi+u\circ\tau$ be the associated normalized observable. Then the definition of the Peierls barrier (see definition~\ref{ManePeierlsBarrier})
implies
\[
h_B(\mathbf{x},\bar{\mathbf{x}}) = h_A(\mathbf{x},\bar{\mathbf{x}}) -u(\bar{\mathbf{x}}) + u(\mathbf{x}),
\quad
\forall\ \mathbf{x},\bar{\mathbf{x}} \in \Sigma.
\]
Since $ h_B(\mathbf{x},\bar{\mathbf{x}}) \geq 0 $, we see that $\mathbf{x}\sim\bar{\mathbf{x}} \Leftrightarrow h_B(\mathbf{x},\bar{\mathbf{x}})=0$ and $h_B(\bar{\mathbf{x}},\mathbf{x})=0$.

To show the transitivity property, it is enough to prove
\[
\mathbf{x}\sim\bar{\mathbf{x}}
\textrm{ and }
\bar{\mathbf{x}}\sim\bar{\bar{\mathbf{x}}}
\Longrightarrow
h_B(\mathbf{x},\bar{\bar{\mathbf{x}}}) = 0.
\]
But proposition \ref{propriedadesbasicas} guarantees
\[
0 \leq h_B(\mathbf{x},\bar{\bar{\mathbf{x}}}) \leq h_B(\mathbf{x},\bar{\mathbf{x}}) +  h_B(\bar{\mathbf{x}},\bar{\bar{\mathbf{x}}})  = 0.
\]
The transitivity property is proved.
\end{proof}

\begin{proposition}\label{proposicao componentes}
The irreducible components are closed and $\sigma$-invariant.
\end{proposition}

\begin{proof}
\noindent{\it Part 1.} Let $\mathbf{x}\in\Omega(A)$. Consider $\{\bar{\mathbf{x}}_\epsilon\}_\epsilon$ a sequence of points of $\Omega(A)$ equivalent to $\mathbf{x}$ and within $\epsilon$ of $\bar{\mathbf{x}}\in\Omega(A)$. Then on the one hand, $h_A(\mathbf{x},\bar{\mathbf{x}}) + h_A(\bar{\mathbf{x}},\mathbf{x}) \geq h_A(\mathbf{x},\mathbf{x})=0$, and on the other hand,
\[
h_A(\mathbf{x},\bar{\mathbf{x}}_\epsilon) + h_A(\bar{\mathbf{x}},\mathbf{x}) \leq h_A(\mathbf{x},\bar{\mathbf{x}}_\epsilon) + h_A(\bar{\mathbf{x}},\bar{\mathbf{x}}_\epsilon) + h_A(\bar{\mathbf{x}}_\epsilon,\mathbf{x}) = h_A(\bar{\mathbf{x}},\bar{\mathbf{x}}_\epsilon).
\]
By continuity of $h_A(\mathbf{x},\cdot)$ and $h_A(\bar{\mathbf{x}},\cdot)$ with respect to the second variable, the previous inequality gives $h_A(\mathbf{x},\bar{\mathbf{x}}) + h_A(\bar{\mathbf{x}},\mathbf{x}) \leq 0$. Therefore $\bar{\mathbf{x}}\sim\mathbf{x}$ and the class containing $\mathbf{x}$ is closed.

\medskip
\noindent{\it Part 2.} Let $\mathbf{x}\in\Omega(A)$. Either $\sigma(\mathbf{x})=\mathbf{x}$ and in an obvious way $\sigma(\mathbf{x})\sim\mathbf{x}$ or $\sigma(\mathbf{x})\not=\mathbf{x}$ and item 3 of proposition~\ref{propriedadesbasicas} shows $\phi_A(\mathbf{x},\sigma(\mathbf{x})) + \phi_A(\sigma(\mathbf{x}),\mathbf{x}) = \phi_A(\mathbf{x},\mathbf{x})=0$. Remember that $h_A(\mathbf{y},\cdot)=\phi_A(\mathbf{y},\cdot)$ whenever $\mathbf{y}\in\Omega(A)$; note that $ \mathbf x $ and $ \sigma(\mathbf x) $ belong to the $\sigma$-invariant set $ \Omega(A) $. Then we get
$ h_A(\mathbf{x},\sigma(\mathbf{x})) + h_A(\sigma(\mathbf{x}),\mathbf{x}) = h_A(\mathbf{x},\mathbf{x})=0 $ and $ \mathbf x $ and $ \sigma(\mathbf x) $ belong to the
same irreducible class.
\end{proof}

We assume from now on that $\Omega(A)$ is equal to a disjoint union of irreducible components, $\Omega(A)=C_1\sqcup\ldots\sqcup C_r$. The following proposition shows that the Peierls barrier normalized by a separating sub-action could play the role of a quantized set of levels of energy.

\begin{proposition}
Let $A$ be a H\"older observable and assume that $\Omega(A)=\sqcup_{i=1}^rC_i$ is equal to a finite union of irreducible components.
\begin{enumerate}
\item If $u$ is a continuous sub-action, then
\[
(\mathbf{x}^i,\mathbf{x}^j)\mapsto h_A(\mathbf{x}^i,\mathbf{x}^j)-u(\mathbf{x}^j)+u(\mathbf{x}^i) \textrm{ is constant on }  C_i\times C_j.
\]
\item If $u$ is a continuous separating sub-action, then
\[
h_A(\mathbf{x}^i,\mathbf{x}^j) > u(\mathbf{x}^j)-u(\mathbf{x}^i),
\quad
\forall\ (\mathbf{x}^i,\mathbf{x}^j) \in C_i\times C_j, \;\; \forall\ i \not= j.
\]
\end{enumerate}
\end{proposition}

\begin{proof}
We first normalize $A$ by taking $B=A-\bar A - u\circ\pi + u\circ\tau$ so that $B\geq 0$ and $\bar B=0$.

\medskip
\noindent{\it Part 1.} Let $(\mathbf{x}^i,\mathbf{x}^j), (\bar{\mathbf{x}}^i,\bar{\mathbf{x}}^j)\in C_i\times C_j$. Then $h_B(\mathbf{x}^i,\bar{\mathbf{x}}^i)=h_B(\mathbf{x}^j,\bar{\mathbf{x}}^j)=0$ and
\[
h_B(\bar{\mathbf{x}}^i,\bar{\mathbf{x}}^j) \leq h_B(\bar{\mathbf{x}}^i,\mathbf{x}^i) + h_B(\mathbf{x}^i,\mathbf{x}^j) + h_B(\mathbf{x}^j,\bar{\mathbf{x}}^j) \leq h_B(\mathbf{x}^i,\mathbf{x}^j).
\]
\noindent Conversely $h_B(\mathbf{x}^i,\mathbf{x}^j) \leq h_B(\bar{\mathbf{x}}^i,\bar{\mathbf{x}}^j)$ and we have proved that $h_B(\cdot,\cdot)$ is constant on $C_i\times C_j$.

\medskip
\noindent{\it Part 2.} Let $\{U_i^\eta\}_{\eta>0}$ be a basis of neighborhoods of $C_i$. Since $\sigma(C_i)\subset C_i$ is disjoint from each $C_j$, $j\not= i$, there exists $\eta>0$ small enough such that $\sigma(U_i^\eta)$ is disjoint from $\cup_{j\not= i}U_j^\eta$. Let $i\not= j$ and $\mathbf{x}\in C_i$, $\bar{\mathbf{x}} \in C_j$. For $\epsilon>0$ sufficiently small, the ball of radius $\epsilon$ centered at $\mathbf{x}$ is included in $U_i^\eta$.
Let $(\mathbf{z}^0 = (\mathbf y^0, \mathbf x^0),\ldots,\mathbf{z}^k = (\mathbf y^k, \mathbf x^k))$ be a path of length $k$ within $\epsilon$ of $\mathbf{x}$
and $\bar{\mathbf{x}}$, more precisely, satisfying $d(\mathbf{x}^0,\mathbf{x})<\epsilon$ and $d(\mathbf{x}^k,\bar{\mathbf{x}})<\epsilon$.
Let $p\geq1$ be  the first time $\sigma^p(\mathbf{x})\not\in U_i^\eta$. Then $\sigma^{p-1}(\mathbf{x})\in U_i^\eta$ and $\sigma^p(\mathbf{x})\in \sigma(U_i^\eta) \setminus U_i^\eta$. By the choice of $\eta$, $\sigma^p(\mathbf{x})\not\in \cup_{j=1}^r U_j^\eta =: \mathcal{U}\supset\Omega(A)$. Since $\Omega(A)=\pi(\mathbb{M}_A(u))$, let $\hat{\mathcal{U}} := \pi^{-1}(\mathcal{U})$, then $\mathbf{z}^p\not\in\hat{\mathcal{U}}$ and
\[
\sum_{l=1}^k B(\mathbf{z}^l) \geq B(\mathbf{z}^p) \geq \min_{\hat{\Sigma}\setminus \hat{\mathcal{U}}} B =: m >0.
\]
\noindent We have proved that $h_B(\mathbf{x},\bar{\mathbf{x}}) \geq m >0$.
\end{proof}

We are now in a position to prove our second result.

\begin{proof}[Proof of theorem \ref{discretestructure}] We fixed once for all $\bar{\mathbf{x}}^i\in C_i$.

\medskip
\noindent{\it Part 1.} We know from theorem \ref{structure} that a continuous calibrated sub-action satisfies
$u(\mathbf{x})=\min_{\bar{\mathbf{x}}\in\Omega(A)} [ u(\bar{\mathbf{x}}) + h_A(\bar{\mathbf{x}},\mathbf{x}) ]$.
If $\bar{\mathbf{x}}\in C_i$, then $\bar{\mathbf{x}}\sim\bar{\mathbf{x}}^i$ and $h_A(\bar{\mathbf{x}}^i,\bar{\mathbf{x}}) + h_A(\bar{\mathbf{x}},\bar{\mathbf{x}}^i)=0$. Then
\begin{align*}
u(\bar{\mathbf{x}}^i) + h_A(\bar{\mathbf{x}}^i,\mathbf{x})
&\leq u(\bar{\mathbf{x}}^i) + h_A(\bar{\mathbf{x}}^i,\bar{\mathbf{x}}) + h_A(\bar{\mathbf{x}},\mathbf{x}) \\
& = u(\bar{\mathbf{x}}^i) - h_A(\bar{\mathbf{x}},\bar{\mathbf{x}}^i) + h_A(\bar{\mathbf{x}},\mathbf{x}) \leq u(\bar{\mathbf{x}}) + h_A(\bar{\mathbf{x}},\mathbf{x}).
\end{align*}
We have proved that $u(\mathbf{x}) = \min_{1\leq i\leq r} [ u(\bar{\mathbf{x}}^i) + h_A(\bar{\mathbf{x}}^i,\mathbf{x})]$. The fact that $(u(\bar{\mathbf{x}}^1),\ldots,u(\bar{\mathbf{x}}^r))\in\mathcal{C}_A(\bar{\mathbf{x}}^1,\ldots,\bar{\mathbf{x}}^r)$ comes from items 1 and 6 of
proposition~\ref{propriedadesbasicas}.

\medskip
\noindent{\it Part 2.} Let $(u_1,\ldots,u_r)\in\mathcal{C}_A(\bar{\mathbf{x}}^1,\ldots,\bar{\mathbf{x}}^r)$ and define $\phi:\Omega(A)\to\mathbb{R}$ by $\phi(\mathbf{x}) := u_i+h_A(\bar{\mathbf{x}}^i,\mathbf{x})$ for all $\mathbf{x}\in C_i$. We notice that $\phi$ is continuous and we show that $\phi(\bar{\mathbf{x}}) - \phi(\mathbf{x}) \leq h_A(\mathbf{x},\bar{\mathbf{x}})$ for all $\mathbf{x},\bar{\mathbf{x}}\in\Omega(A)$. Indeed, if $\mathbf{x}\in C_i$ and $\bar{\mathbf{x}}\in C_j$, then
\begin{align*}
\phi(\bar{\mathbf{x}}) - \phi(\mathbf{x}) &= (u_j-u_i) + h_A(\bar{\mathbf{x}}^j,\bar{\mathbf{x}}) - h_A(\bar{\mathbf{x}}^i,\mathbf{x}) \\
&\leq h_A(\bar{\mathbf{x}}^i,\bar{\mathbf{x}}^j) + h_A(\bar{\mathbf{x}}^j,\bar{\mathbf{x}}) - h_A(\bar{\mathbf{x}}^i,\mathbf{x}) \\
&= h_A(\bar{\mathbf{x}}^i,\bar{\mathbf{x}}^j) - h_A(\bar{\mathbf{x}},\bar{\mathbf{x}}^j) - h_A(\bar{\mathbf{x}}^i,\mathbf{x}) \\
&\leq h_A(\bar{\mathbf{x}}^i,\bar{\mathbf{x}}) - h_A(\bar{\mathbf{x}}^i,\mathbf{x}) \leq h_A(\mathbf{x},\bar{\mathbf{x}}).
\end{align*}
\noindent (The last but one inequality uses item 1 of proposition \ref{propriedadesbasicas} and the fact that $h_A(\bar{\mathbf{x}}^i,\cdot)$ is a sub-action.)
By theorem~\ref{structure}, we know that the function
$ u(\mathbf{x}) := \min_{\bar{\mathbf{x}}\in\Omega(A)} [ \phi(\bar{\mathbf{x}}) + h_A(\bar{\mathbf{x}},\mathbf{x}) ]$ is
a continuous calibrated sub-action which extends $\phi$ on $\Omega(A)$. In particular, $u(\bar{\mathbf{x}}^i)=\phi(\bar{\mathbf{x}}^i)=u_i$ and, thanks to part 1, $u$ coincides with $\min_{1\leq i\leq r} [ u_i + h_A(\bar{\mathbf{x}}^i,\cdot) ]$.

\medskip
\noindent{\it Part 3.} Let $i_0 \in \{1, \ldots, r\} $. If $(u_1,\ldots,u_r)$ satisfies $u_i = u_{i_0} + h_A(\bar{\mathbf{x}}^{i_0},\bar{\mathbf{x}}^i)$,
then $i_0$ is unique. Otherwise there would exist $i_1 \not= i_0$ such that $u_i = u_{i_1} + h_A(\bar{\mathbf{x}}^{i_1},\bar{\mathbf{x}}^i)$. Thus
\[
u_{i_1} = u_{i_0} + h_A(\bar{\mathbf{x}}^{i_0},\bar{\mathbf{x}}^{i_1})
\quad\textrm{and}\quad
u_{i_0} = u_{i_1} + h_A(\bar{\mathbf{x}}^{i_1},\bar{\mathbf{x}}^{i_0}).
\]
We would obtain $h_A(\bar{\mathbf{x}}^{i_0},\bar{\mathbf{x}}^{i_1}) + h_A(\bar{\mathbf{x}}^{i_1},\bar{\mathbf{x}}^{i_0}) = 0$ contradicting $\bar{\mathbf{x}}^{i_0} \not\sim \bar{\mathbf{x}}^{i_1}$. The fact that $(u_1,\ldots,u_r)\in\mathcal{C}_A(\bar{\mathbf{x}}^1,\ldots,\bar{\mathbf{x}}^r)$ comes from
\[
u_j-u_i = h_A(\bar{\mathbf{x}}^{i_0},\bar{\mathbf{x}}^{j}) - h_A(\bar{\mathbf{x}}^{i_0},\bar{\mathbf{x}}^{i}) \leq h_A(\bar{\mathbf{x}}^{i},\bar{\mathbf{x}}^{j}).
\]
\noindent The end of part 3 follows since $u(\mathbf{x}) := u_{i_0} + h_A(\bar{\mathbf{x}}^{i_0},\mathbf{x})$ already defines a calibrated sub-action satisfying $u(\bar{\mathbf{x}}^i)=u_{i}$ for all $i$.
\end{proof}

The proof of application \ref{application} is elementary.

\begin{proof}[Proof of application \ref{application}] Define $B := A-v\circ\pi + v\circ\tau-\bar A$, then the null function is a sub-action of $B$ and $v-u$ is a sub-action calibrated to $B$. Moreover, $h_B(\mathbf{x},\bar{\mathbf{x}}) = h_A(\mathbf{x},\bar{\mathbf{x}}) - v(\bar{\mathbf{x}}) + v(\mathbf{x})$ and $\Omega(A)=\Omega(B)$. It is therefore enough to assume $A$ normalized ($A\geq0$ and $\bar A=0$) and $v=0$.

\medskip
\noindent{\it Part 1.} If $\mathbf{x}\sim\bar{\mathbf{x}}$ are two points of $\Omega(A)$, then $h_A(\mathbf{x},\bar{\mathbf{x}})=0$ and $h_A(\bar{\mathbf{x}},\mathbf{x})=0$. Thanks to items 1 and 6 of proposition \ref{propriedadesbasicas}, we obtain $u(\mathbf{x})=u(\bar{\mathbf{x}})$. If $\mathbf{x}$ is any point of $\Sigma$, by the calibration of $u$, one can construct an inverse path $\{\mathbf{z}^{-i}\}_{i\geq0}$ of $\hat{\Sigma}$,
with $\pi(\mathbf{z}^0)=\mathbf{x}$,
such that $u(\mathbf{x}^{-i})-u(\mathbf{x}^{-i-1})=A(\mathbf{z}^{-i})$, $\mathbf{x}^{-i}=\pi(\mathbf{z}^{-i})$, for all $i$. Let $\bar{\mathbf{x}}$ be an accumulation point of $\{\mathbf{x}^{-i}\}_{i\geq0}$. Then $\bar{\mathbf{x}}\in\Omega(A)$ and, since $A\geq0$, the sequence $\{u(\mathbf{x}^{-i})\}_{i\geq0}$ is decreasing.
In particular, $u(\mathbf{x})\geq u(\bar{\mathbf{x}}) $ establishes $ \min_\Sigma u = \min_{\Omega(A)} u $.

\medskip
\noindent{\it Part 2.} Let $u_i$ be the value of $u$ on $C_i$. Assume we have ordered these values as $u_{i_1}\leq u_{i_2} \leq \ldots \leq u_{i_r}$. Let $\bar{\mathbf{x}}^i\in C_i$ fixed. It suffices to prove $u(\bar{\mathbf{x}}^{i_k}) = u(\bar{\mathbf{x}}^{i_1}) + h_A(\bar{\mathbf{x}}^{i_1},\bar{\mathbf{x}}^{i_k})$ for all $k=1,\ldots,r$. It is true for $k=1$. Since $C_{i_{k+1}}$ is not a minimum local of $u$, one can find a sequence of points $\{\mathbf{x}_\epsilon\}_{\epsilon>0}$ within $\epsilon$ of $C_{i_{k+1}}$ such that $u(\mathbf{x}_\epsilon) < u(\bar{\mathbf{x}}^{i_{k+1}})$. From part 1 of theorem \ref{discretestructure}, there exists an index $j$ such that $u(\mathbf{x_\epsilon}) = u(\bar{\mathbf{x}}^{j}) + h_A(\bar{\mathbf{x}}^j,\mathbf{x}_\epsilon)$. Since $h_A\geq0$, $u_j=u(\bar{\mathbf{x}}^{j}) \leq u(\mathbf{x_\epsilon}) < u_{i_{k+1}}$. So $j$ has to be one of indexes $i_1,\ldots,i_k$. By induction, $u(\bar{\mathbf{x}}^{j}) = u(\bar{\mathbf{x}}^{i_1}) + h_A(\bar{\mathbf{x}}^{i_1},\bar{\mathbf{x}}^{j})$ and
\[
u(\mathbf{x_\epsilon}) = u(\bar{\mathbf{x}}^{i_1}) + h_A(\bar{\mathbf{x}}^{i_1},\bar{\mathbf{x}}^{j}) + h_A(\bar{\mathbf{x}}^j,\mathbf{x}_\epsilon).
\]
On the one hand,
$h_A(\bar{\mathbf{x}}^{i_1},\bar{\mathbf{x}}^{j}) + h_A(\bar{\mathbf{x}}^j,\mathbf{x}_\epsilon) \geq h_A(\bar{\mathbf{x}}^{i_1},\mathbf{x}_\epsilon)$ implies
\[
u(\mathbf{x_\epsilon}) \geq u(\bar{\mathbf{x}}^{i_1}) + h_A(\bar{\mathbf{x}}^{i_1},\mathbf{x}_\epsilon).
\]
On the other hand, as $u$ is a sub-action, we obtain the reverse inequality and finally
\[
u(\mathbf{x_\epsilon}) = u(\bar{\mathbf{x}}^{i_1}) + h_A(\bar{\mathbf{x}}^{i_1},\mathbf{x}_\epsilon).
\]
Letting $\epsilon$ go to 0, $\mathbf{x}_\epsilon$ accumulates to $C_{i_{k+1}}$ and
\[
u(\bar{\mathbf{x}}^{i_{k+1}}) = u(\bar{\mathbf{x}}^{i_1}) + h_A(\bar{\mathbf{x}}^{i_1},\bar{\mathbf{x}}^{i_{k+1}}).
\]
\end{proof}
\end{section}

\end{document}